\documentclass[a4paper,10pt]{article}
\usepackage{amsfonts}
\usepackage{amsmath}
\usepackage{amssymb}
\usepackage{theorem}
\usepackage{amscd}
\usepackage{pstricks}
\usepackage[english]{babel}
\newtheorem{thm}{Theorem}[section]

\newtheorem{df}[thm]{Definition}

\newtheorem{prop}[thm]{Proposition}

\newcommand{\Irr}[1]{\operatorname{Irr}(#1)}\newcommand{\res}[1]{\operatorname{res}(#1)}
\newcommand{\rad}[1]{\operatorname{rad}(#1)}

\def\[{[\![}
\def\]{]\!]}
\title{An algorithm for the   computation of the decomposition matrices  for Ariki-Koike algebras}
\author{Nicolas JACON \footnote{Address: {\sc Institut Girard Desargues, bat. Jean Braconnier, Universit\'e
Lyon 1, 21 av Claude Bernard, F--69622 Villeurbanne cedex, France.} E-mail address: jacon@igd.univ-lyon1.fr}}
\date{}
\begin{document}
\maketitle
{\abstract We give a purely combinatorial algorithm for the computation of the decomposition matrices for Ariki-Koike algebras when the parameters are powers of the same root of unity.}

\section{Introduction}
Ariki-Koike algebras have been independantly introduced by Ariki and Koike in \cite{AriKoi} and by Brou\'e and Malle in \cite{BrouMAl}. According to a conjecture of   Brou\'e and Malle, this kind of algebras should play a role in the decomposition of the induced cuspidal representations of the finite groups of Lie type. 
Let  $R$ be a commutative ring, let $d\in{\mathbb{N}_{>0}}$,  $n\in{\mathbb{N}}$ and let  $v$, $u_0$, $u_1$,..., $u_{d-1}$ be  $d+1$ parameters in $R$. We consider the Ariki-Koike algebra $\mathcal{H}_{R,n}:=\mathcal{H}_{R,n}(v;u_0,...,u_{d-1})$ of type $G(d,1,n)$ over $R$. This is the unital associative $R$-algebra defined by:
\begin{itemize}
\item generators: $T_0$, $T_1$,..., $T_{n-1}$,
\item relations symbolized by the following diagram:
\\
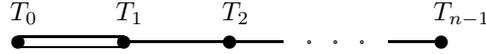
\begin{center}
\begin{picture}(240,20)
\put( 50,10){\circle*{5}}
\put( 47,18){$T_0$}
\put( 50,8){\line(1,0){40}}
\put( 50,12){\line(1,0){40}}
\put( 90,10){\circle*{5}}
\put( 87,18){$T_1$}
\put( 90,10){\line(1,0){40}}
\put(130,10){\circle*{5}}
\put(127,18){$T_2$}
\put(130,10){\line(1,0){20}}
\put(160,10){\circle{1}}
\put(170,10){\circle{1}}
\put(180,10){\circle{1}}
\put(190,10){\line(1,0){20}}
\put(210,10){\circle*{5}}
\put(207,18){$T_{n-1}$}
\end{picture}
\end{center}
and the following ones:
 \begin{align*}
                &(T_0-u_0)(T_0-u_1)...(T_0-u_{d-1})  =  0,\\
                             &(T_i-v)(T_i+1)  =  0\ (i\geq{1}). \end{align*} 
\end{itemize}

Assume that $R$ is a field of characteristic $0$. Let $\Pi^d_n$ be the set of $d$-partitions of rank $n$ that is to say the set of  $d$-tuples of partitions  $\underline{\lambda}=(\lambda^{(0)},...,\lambda^{(d-1)})$ such that $|\displaystyle{\lambda^{(0)}|+...+|\lambda^{(d-1)}|=n}$.
For each  $\underline{\lambda}\in{\Pi^d_n} $, Dipper, James and Mathas (\cite{DJM}) have  defined a right $\mathcal{H}_{R,n}$-module $S^{\underline{\lambda}}_R$ which is called a Specht module\footnote{ Here, we use the definition of the classical Specht modules. Note that the  results in  \cite{DJM} are in fact given in terms of dual Specht modules. The passage from classical Specht modules to their duals is provided by the map $(\lambda^{(0)},\lambda^{(1)},...,\lambda^{(d-1)})\mapsto{(\lambda^{(d-1)'},\lambda^{(d-2)'},...,\lambda^{(0)'})}$ where, for $i=0,...,d-1$,  $\lambda^{(i)'}$ is denoting the conjugate partition.}. For each Specht module  $S^{\underline{\lambda}}_R$, they have  attached a natural bilinear form and a radical $\rad{S^{\underline{\lambda}}_R}$ such that the non zero $D^{\underline{\lambda}}_R:=S^{\underline{\lambda}}_R/\rad{S^{\underline{\lambda}}_R}$ form a complete set of non isomorphic irreducible modules. Let  $\Phi^d_n:=\{\underline{\mu}\in{\Pi^d_n}\ |\ D_R^{\underline{\mu}}\neq{0}\}$.

Let $R_0{(\mathcal{H}_{R,n})}$ be the Grothendieck group of finitely generated  $\mathcal{H}_{R,n}$-modules. This is generated by the set of simple $\mathcal{H}_{R,n}$-modules. Thus, for each $\underline{\lambda}\in{\Pi_n^d}$ and $\underline{\mu}\in{\Phi_n^d}$, there exist numbers $d_{\underline{\lambda},\underline{\mu}}$ which are called the decomposition numbers such that:
$$[S^{\underline{\lambda}}_R]=\sum_{\underline{\mu}\in{\Phi_n^d}}d_{\underline{\lambda},\underline{\mu}}[D^{\underline{\mu}}_R].$$
The matrix $(d_{\underline{\lambda},\underline{\mu}})_{\underline{\lambda}\in{\Pi^d_n},\underline{\mu}\in{\Phi^d_n}}$ is called the decomposition matrix of $\mathcal{H}_{R,n}$.

One of the main problems in the representation theory of Ariki-Koike algebras is the determination of the decomposition matrix. When $\mathcal{H}_{R,n}$ is semi-simple, the decomposition matrix is just the identity. When $\mathcal{H}_{R,n}$ is not semi-simple, by using results of Dipper and Mathas, the determination of the decomposition matrix is deduced from the case where all the parameters are powers of the same number $\eta$ (see \cite{DipMat}). Here, we assume that $\eta$ is a primitive $e^{\textrm{th}}$-root of unity.

 When $d=1$, Lascoux, Leclerc and Thibon \cite{LLT} have presented a fast algorithm for the computation of the canonical basis elements of a certain integrable $\mathcal{U}_q(\widehat{sl_e})$-module $\overline{\mathcal{M}}$. Moreover, they conjectured that the problem of computing  the decomposition matrix of $\mathcal{H}_{R,n}$ can be translated to that of computing the canonical basis of  $\overline{\mathcal{M}}$. 
This conjecture has been proved and generalized for all $d\in{\mathbb{N}_{>0}}$ by Ariki in \cite{Ari4}.  Unfortunately, the generalization doesn't give an analogue of the LLT algorithm for $d>1$. In this case,  Uglov \cite{uglov} has given an algorithm but it computes the canonical basis for a larger space which contains $\overline{\mathcal{M}}$ as a submodule. It might be interesting to obtain a generalization of the LLT algorithm for all $d>1$.

 In \cite{these} and \cite{cyclo}, extending the results  developed in \cite{geckrouq} for $d=1$ and $d=2$ by using an ordering of Specht modules by Lusztig $a$-function, we showed that there exists a ``canonical basic set'' $\mathcal{B}$ of Specht modules in bijection with $\Irr{\mathcal{H}_{R,n}}$ and that this set is parametrized by some FLOTW $d$-partitions defined by Foda et al. in \cite{FLOTW}. As a consequence, this result gives a purely combinatorial triangular  algorithm for the computation of the decomposition matrix for Ariki-Koike algebras which generalizes the LLT algorithm.

The aim of this paper is to present this  algorithm. In the first part, we give  the  definitions and theorems used in the algorithm. Then, we give the different steps of the algorithm. 
\section{Ariki's theorem and canonical basic set}

\indent Let $e$ and $d$ be two positive integers and  let $\displaystyle{\eta_e:=\textrm{exp}(\frac{2i\pi}{e})}$ and  $\displaystyle{\eta_d:=\textrm{exp}(\frac{2i\pi}{d})}$. We consider the Ariki-Koike algebra $\mathcal{H}_{R,n}$ over $R:=\mathbb{Q}[\eta_d](\eta_e)$ with the following choice of parameters:
$$v=\eta_e, \qquad{u_j=\eta_e^{v_j},}\qquad{\textrm{for}\ j=0,...,d-1},$$
where $0\leq{v_0}\leq ... \leq{v_{d-1}}<e$. In this section, we briefly summarize the results of Ariki which  give an interpretation of the decomposition matrix in terms of the canonical basis of a certain $\mathcal{U}_q(\widehat{\textrm{sl}_e})$-module. For more details, we refer to  \cite{Arilivre} and to \cite{mathas}. Next,  we recall the results shown in   \cite[chapter 2]{these} and  in \cite{cyclo}.
 \\
a) We first explain the Ariki's theorem. To do this, we need some combinatorial definitions.  Let   $\underline\lambda={(\lambda^{(0)} ,...,\lambda^{(d-1)})}$ be a $d$-partition of rank  $n$. The diagram of  $\underline{\lambda}$ is the following set:
$$[\underline{\lambda}]=\left\{ (a,b,c)\ |\ 0\leq{c}\leq{d-1},\ 1\leq{b}\leq{\lambda_a^{(c)}}\right\}.$$
The elements of this diagram are called   the nodes of  $\underline{\lambda}$.
Let  $\gamma=(a,b,c)$ be a node of  $\underline{\lambda}$. The residue of  $\gamma$  associated to the set  $\{e;{v_0},...,{v_{d-1}}\}$ is the element of $\mathbb{Z}/e\mathbb{Z}$ defined by:
$$\textrm{res}{(\gamma)}=(b-a+v_{c})(\textrm{mod}\ e).$$
If $\gamma $ is a node with residue $i$, we say that  $\gamma$ is an  $i$-node.
Let  $\underline{\lambda}$ and  $\underline{\mu}$ be two $d$-partitions of rank  $n$ and $n+1$ such that  $[\underline{\lambda}]\subset{[\underline{\mu}]}$. There exists a node $\gamma$ such that  $[\underline{\mu}]=[\underline{\lambda}]\cup{\{\gamma\}}$. Then, we denote $[\underline{\mu}]/[\underline{\lambda}]=\gamma$ and  if $\textrm{res}{(\gamma)}=i$, we say that  $\gamma$ is an addable $i$-node  for   $\underline{\lambda}$ and a removable  $i$-node for  $\underline{\mu}$. Now, we introduce an order on the set of nodes of a $d$-partition.
 We say that  $\gamma=(a,b,c)$ is above  $\gamma'=(a',b',c')$ if:
$$b-a+v_c<b'-a'+v_{c'}\ \textrm{or } \textrm{if}\ b-a+v_c=b'-a'+v_{c'}\textrm{ and }c>c'.$$
Let $\underline{\lambda}$ and $\underline{\mu}$ be two $d$-partitions of rank $n$ and $n+1$ such that  there exists an $i$-node $\gamma$ such that  $[\underline{\mu}]=[\underline{\lambda}]\cup{\{\gamma\}}$. We define the following numbers:
\begin{align*}
\overline{N}_i^{a}{(\underline{\lambda},\underline{\mu})}=&   \sharp\{ \textrm{addable }\ i-\textrm{nodes of } \underline{\lambda}\ \textrm{ above } \gamma\} \\
 & -\sharp\{ \textrm{removable }\ i-\textrm{nodes of } \underline{\mu}\ \textrm{ above } \gamma\},\\
\overline{N}_i^{b}{(\underline{\lambda},\underline{\mu})}= &   \sharp\{ \textrm{addable } i-\textrm{nodes of } \underline{\lambda}\ \textrm{ below } \gamma\}\\
&   -\sharp\{ \textrm{removable } i-\textrm{nodes of } \underline{\mu}\ \textrm{ below } \gamma\},\\
\overline{N}_{i}{(\underline{\lambda})} =& \sharp\{ \textrm{addable } i-\textrm{nodes of } \underline{\lambda}\}\\
& -\sharp\{ \textrm{removable } i-\textrm{nodes of } \underline{\lambda}\},\\
\overline{N}_{\mathfrak{d}}{(\underline{\lambda})} =& \sharp\{ 0-\textrm{nodes of } \underline{\lambda}\}.
\end{align*}
b)   Now, let $\mathfrak{h}$ be the free $\mathbb{Z}$-module with basis $\{h_i,\mathfrak{d}\ |\ 0\leq i<e\}$ as in \cite[section 2.B]{cyclo}, let $q$ be an indeterminate and let $\mathcal{U}_q:=\mathcal{U}_q(\widehat{\textrm{sl}_e})$ be the quantum group of type $A^{(1)}_{e-1}$. This is a unital associative algebra over  $\mathbb{C}(q)$ which is generated by elements $\{e_i,f_i\ |\ i\in{\{0,...,e-1\}}\}$ and $\{k_h\ |\ h\in{\mathfrak{h}}\}$ (see \cite[Definition 3.16]{Arilivre} for the relations). Let $\mathcal{A}=\mathbb{Z}[q,q^{-1}]$. We consider the Kostant-Lusztig form of $\mathcal{U}_q$  which is denoted by  $\mathcal{U}_{\mathcal{A}}$: this is a  $\mathcal{A}$-subalgebra of $\mathcal{U}_q$ generated by the  divided powers  $e_i^{(r)}$, $f_j^{(r)}$  for  $0\leq{i,j}<e$, $r\in{\mathbb{N}}$ and by  $k_{h_i}$,  $k_{\mathfrak{d}}$,  $k_{h_i}^{-1}$,  $k^{-1}_{\mathfrak{d}}$ for  $0\leq{i}<e$. Now, if $S$ is a ring and $u$   an invertible element in $S$, we can form the specialized algebra $\mathcal{U}_{S,u}:=S\otimes_{\mathcal{A}}\mathcal{U}_{\mathcal{A}}$  by specializing the indeterminate $q$ io $u\in{S}$.

 For $n\in{\mathbb{N}}$, let $\mathcal{F}_n:=\{\underline{\lambda}\ |\  \underline{\lambda}\in{\Pi_d^n}\}$ and let  $\mathcal{F}:=\bigoplus_{n\in{\mathbb{N}}}  \mathcal{F}_n$.  $\mathcal{F}$ is called the Fock space.
Then, the following theorem shows that we have a  $\mathcal{U}_q$-module structure on  $\mathcal{F}$.
\begin{thm}\label{jmmo} (Jimbo, Misra, Miwa, Okado \cite{JMMO})   $\mathcal{F}$ is a  $\mathcal{U}_q$-module with action:
$$e_{i}\underline{\lambda}=\sum_{\res{[\underline{\lambda}]/[\underline{\mu}]}=i}{q^{-\overline{N}_i^{a}{(\underline{\mu},\underline{\lambda})}}\underline{\mu}},\qquad{f_{i}\underline{\lambda}=\sum_{\res{[\underline{\mu}]/[\underline{\lambda}]}=i}{q^{\overline{N}_i^{b}{(\underline{\lambda},\underline{\mu})}}\underline{\mu}}},$$
$$k_{h_i}\lambda=q^{\overline{N}_i{(\underline{\lambda})}}\underline{\lambda},\qquad{k_{\mathfrak{d}}\underline{\lambda}=q^{-\overline{N}_{\mathfrak{d}}{(\underline{\lambda})}}\underline{\lambda}},$$
 where $0\leq{i}\leq{n-1}$.
\end{thm}
Note that this action is distinct from the action used by Ariki and Mathas for example in \cite{AriMat1}.
Let $\overline{\mathcal{M}}$ be the $\mathcal{U}_q$-submodule of   $\mathcal{F}$ generated by the empty $d$-partition.  This is an integrable  highest weight module. Thus,  we can use the canonical basis theory  to obtain a basis for  $\overline{\mathcal{M}}_{\mathcal{A}}$, the  $\mathcal{U}_{\mathcal{A}}$-module generated by the empty $d$-partition. In particular, the canonical basis elements are known to be indexed by the vertices of some ``crystal graph''. In \cite{FLOTW}, Foda, Leclerc, Okado, Thibon and Welsh have shown that the vertices of the crystal graph of  $\overline{\mathcal{M}}$ are labeled by the following $d$-partitions:
\begin{df} (Foda, Leclerc, Okado, Thibon, Welsh \cite{FLOTW})
 We say that  $\underline\lambda={(\lambda^{(0)} ,...,\lambda^{(d-1)})}$
is a FLOTW  $d$-partition associated to  the set  ${\{e;{v_0},...,v_{d-1}\}}$ if and only if:
\begin{enumerate}
\item for all $0\leq{j}\leq{d-2}$ and $i=1,2,...$, we have:
\begin{align*}
&\lambda_i^{(j)}\geq{\lambda^{(j+1)}_{i+v_{j+1}-v_j}},\\
&\lambda^{(d-1)}_i\geq{\lambda^{(0)}_{i+e+v_0-v_{d-1}}};
\end{align*}
\item  for all  $k>0$, among the residues appearing at the right ends of the length $k$ rows of   $\underline\lambda$, at least one element of  $\{0,1,...,e-1\}$ does not occur.
\end{enumerate}
We denote by $\Lambda^{1,n}_{\{e;{v_0},...,{v_{d-1}}\}}$ the set of FLOTW $d$-partitions with rank $n$ associated to  ${\{e;{v_0},...,{v_{d-1}}\}}$. If there is no ambiguity concerning  $\{e;{v_0},...,{v_{d-1}}\}$, we denote it by $\Lambda^{1,n}$.
\end{df}
Now the canonical basis of $\overline{\mathcal{M}}$ is defined by using the following theorem:
\begin{thm}\label{basedef} (Kashiwara-Lusztig, see \cite[chapter 9]{Arilivre}) Define the bar involution to be the $\mathbb{Z}$- linear ring automorphism of $\mathcal{U}_{\mathcal{A}}$ determined for $i=0,...,e-1$ and $h\in{\mathfrak{h}}$ by:
$$\overline{q}:=q^{-1},\qquad{\overline{k_h}=k_{-h}},\qquad{\overline{e_i}:=e_i,}\qquad{\overline{f_i}:=f_i}.$$
We extend it to  $\overline{\mathcal{M}_{\mathcal{A}}}$ by setting $\overline{u .\underline{\emptyset}}:=\overline{u}.\underline{\emptyset}$ for all $u\in{\mathcal{U}_{\mathcal{A}}}$.
Then, for each  $\underline{\mu}\in{\Lambda^{1,n}}$,  there exists a unique element $G({\underline{\mu}})$ in  $\overline{\mathcal{M}}_{\mathcal{A}}$ such that:
\begin{itemize}
\item $\overline{ G({\underline{\mu}})}=G({\underline{\mu}}),$
\item $G({\underline{\mu}})=\underline{\mu} \ (\textrm{mod}\  q).$
\end{itemize}
The set  $\{G({\underline{\mu}})\ |\ \underline{\mu}\in{\Lambda^{1,n}}\}$ is a basis  of  $\overline{\mathcal{M}}_{\mathcal{A}}$ which is uniquely determined by the above conditions. It is called the canonical basis of  $\overline{\mathcal{M}}$.
\end{thm}
Now the following result of Ariki which were conjectured by Lascoux, Leclerc and Thibon in \cite{LLT} for $l=1$ gives an interpretation of the decomposition matrix of $\mathcal{H}_{R,n}$ in terms of the canonical basis of  $\overline{\mathcal{M}}$.
\begin{thm}(Ariki \cite{Ari4})  Let $\underline{\mu}\in{\Lambda^{1,n}}$,  there exist polynomials $b_{\underline{\lambda},\underline{\mu}}(q)\in{q\mathbb{Z}[q]}$ such that:
$$G({\underline{\mu}})=\sum_{\underline{\lambda}\in{\Pi_d^n}}b_{\underline{\lambda},\underline{\mu}}(q) \underline{\lambda}.$$
Then,  there exists a unique bijection $j:\Lambda^{1,n}\to{\Phi_n^d}$ 
 such that  $b_{\underline{\lambda},\underline{\mu}}(1)=d_{\underline{\lambda},{j(\underline{\mu})}}$ for all $\underline{\lambda}\in{\Pi_n^d}$ where $(d_{\underline{\lambda},\underline{\nu}})_{\underline{\lambda}\in{\Pi^d_n},\underline{\nu}\in{\Phi^d_n}}$ is the decomposition matrix of   $\mathcal{H}_{R,n}$.
\end{thm}
Thus, the elements of the canonical basis evaluated at $q=1$ correspond to the columns of the decomposition matrix of $\mathcal{H}_{R,n}$ that is to say  the  indecomposable  projective $\mathcal{H}_{R,n}$-modules.\\ 
c) The aim of the work presented in \cite{cyclo} was to study the indecomposable projective $\mathcal{H}_{R,n}$-modules. 
The main result  gives an interpretation of the decomposition matrix   of   $\mathcal{H}_{R,n}$ in terms of Lusztig $a$-function. In particular, extending results  of Geck and Rouquier (see \cite{geckrouq}), we proved that there exists a so called ``canonical basic set'' of Specht modules which is in bijection with the set of simple   $\mathcal{H}_{R,n}$-modules. 

The first step is to define an ``$a$-value'' on each $d$-partition. To do this, we consider a semi-simple  Ariki-Koike algebra of type $G(d,1,n)$ with a special choice of parameters and we define an $a$-value on the simple   modules (which are parametrized by the $d$-partitions) using the characterization of the Schur elements which have been obtained by Geck, Iancu and Malle in \cite{GIM}. This leads to  the following definition:
\begin{df}\label{afonction}  Let $\underline{\lambda}:=(\lambda^{(0)},\lambda^{(1)},...,\lambda^{(d)})\in{\Lambda^{1,n}}$ where for $i=0,...,d-1$ we have $\lambda^{(i)}:=(\lambda^{(i)}_1,...,\lambda^{(i)}_{h^{(i)}})$. We assume that the rank of  $\underline{\lambda}$ is $n$.   For $i=0,...,d-1$ and $p=1,...,n$, we define the following rational numbers:
\begin{align*}
 m^{(i)}&:=v_i-\frac{ie}{d}+e,\\
B'^{(i)}_p&:=\lambda^{(i)}_p-p+n+m^{(i)},\end{align*}
where we use the convention that $\lambda^{(i)}_p:=0$ if $p>h^{(i)}$. For $i=0,...,d-1$, let $B'^{(i)}=(B'^{(i)}_1,...,B'^{(i)}_n)$.   Then, we define:
$$a_1(\underline{\lambda}):=\sum_{{0\leq{i}\leq{j}<d}\atop{{(a,b)\in{B'^{(i)}\times{B'^{(j)}}}}\atop{a>b\ \textrm{if}\ i=j}}}{\min{\{a,b\}}}   -\sum_{{0\leq{i,j}<d}\atop{{a\in{B'^{(i)}}}\atop{1\leq{k}\leq{a}}}}{\min{\{k,m^{(j)}\}}}.$$
\end{df}
 Now, the $a$-value associated to  $S^{\underline{\lambda}}_R$ is the rational number $a(\underline{\lambda}):=a_1(\underline{\lambda})+f(n)$ where $f(n)$ is a rational number which only depends on the parameters $\{e;v_0,...,v_{d-1}\}$ and on $n$ (the expression of $f$ is given in \cite{cyclo}).

Next, we associate to  each  $\underline{\lambda}\in{\Lambda^{1,n}}$ a sequence of residues which will have ``nice'' properties with respect to the $a$-value:
\begin{prop}\label{asuite} (\cite[Definition 4.4]{cyclo}) Let $\underline{\lambda}\in{\Lambda^{1,n}}$ and let:
$$l_{\textrm{max}}:=\operatorname{max}\{\lambda^{(0)}_1,...,\lambda^{(l-1)}_1\}.$$
Then,  there exists a removable node $\xi_1$ with residue  $k$ on a part  $\lambda_{j_1}^{(i_1)}$ with  length $l_{\textrm{max}}$, such that there doesn't exist a $k-1$-node at the right end of a part with length $l_{\textrm{max}}$ (the existence of such a node is proved in \cite[Lemma 4.2]{cyclo}). 

Let $\gamma_1$, $\gamma_2$,..., $\gamma_r$ be the  $k-1$-nodes at the right ends of parts $\displaystyle {\lambda_{p_1}^{(l_1)}\geq{\lambda_{p_2}^{(l_2)}\geq{...}}}$$\geq{\lambda_{p_r}^{(l_r)}}$. 
 Let   $\xi_1$, $\xi_2$,..., $\xi_{s}$ be the  removable $k$-nodes  of   $\underline{\lambda}$ on  parts  $\displaystyle{\lambda_{j_1}^{(i_1)}\geq{\lambda_{j_2}^{(i_2)}\geq{...}\geq{\lambda_{j_s}^{(i_s)}}}}$   such that:
$$\lambda_{j_s}^{(i_s)}>\lambda_{p_1}^{(l_1)}.$$

We remove the nodes  $\xi_1$, $\xi_2$,..., $\xi_s$ from $\underline{\lambda}$. Let  $\underline{\lambda}'$ be the  resulting $d$-partition. Then,   $\underline{\lambda}'\in{\Lambda}^{1,n-s}$ and we define recursively the  $a$-sequence  of residues of  $\underline{\lambda}$ by:
$$a\textrm{-sequence}(\underline{\lambda})=a\textrm{-sequence}(\underline{\lambda}'),\underbrace{k,...,k}_{s}.$$
\end{prop}
\textbf{Example}:\\
Let $e=4$, $d=3$, $v_0=0$, $v_1=2$ and $v_2=3$. We consider the  $3$-partition $\underline{\lambda}=(1,3.1,2.1.1)$ with the following diagram:
$$  \left( \ \begin{tabular}{|c|}
                         \hline
                           0    \\
                          \hline
                             \end{tabular}\  ,\
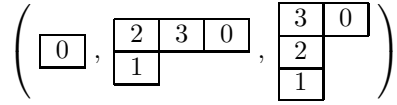
                       \begin{tabular}{|c|c|c|}
                         \hline
                           2  & 3  & 0 \\
                          \hline
                          1     \\
                          \cline{1-1}
                       
                             \end{tabular}\ ,\ 
                        \begin{tabular}{|c|c|}
                         \hline
                           3  & 0 \\
                          \hline
                          2     \\
                          \cline{1-1}
                          1 \\
                           \cline{1-1}
                             \end{tabular}\  
                                                 \right)$$
\\
$\underline{\lambda}$ is a FLOTW $3$-partition.

We search the $a$-sequence of $\underline{\lambda}$: we have to find  $k\in{\{0,1,2,3\}}$, $s\in{\mathbb{N}}$ and a $2$-partition $\underline{\lambda}'$ such that:
$$a\textrm{-sequence}(\underline{\lambda})=a\textrm{-sequence}(\underline{\lambda}'),\underbrace{k,...,k}_{s}.$$
The part with maximal length is  the  part with length  $3$ and the residue of the associated removable node is $0$. We remark that there are two others removable  $0$-nodes on parts with length $1$ and $2$. Since there is no  node with residue $0-1\equiv 3\ (\textrm{mod}\ e)$ at the right ends of the parts of  $\underline{\lambda}$, we must remove these three $0$-nodes. Thus, we have to take   $k=0$, $s=3$ and   $\underline{\lambda}'=(\emptyset,2.1,1.1.1)$, hence:
$$a\textrm{-sequence}(\underline{\lambda})=a\textrm{-sequence}(\emptyset,2.1,1.1.1),0,0,0.$$
Observe  that the  $3$-partition $(\emptyset,2.1,1.1.1)$ is a FLOTW  $3$-partition.

Now, the residue of the removable node on the part with maximal length is $3$. Thus, we obtain:
$$a\textrm{-sequence}(\underline{\lambda})=a\textrm{-sequence}(\emptyset,1.1,1.1.1),3,0,0,0.$$
Repeating the same procedure, we finally  obtain:
$$a\textrm{-sequence}(\underline{\lambda})=3,2,2,1,1,3,0,0,0.$$
\\
\begin{prop}\label{decompo} (\cite[Proposition 4.14]{cyclo}) Let $n\in{\mathbb{N}}$, let $\underline{\lambda}\in{\Lambda}^{1,n}$  and let $a\textrm{-sequence}(\underline{\lambda})=\underbrace{i_1,...,i_1}_{a_1},\underbrace{i_{2},...,i_{2}}_{a_{2}},...,\underbrace{i_s,...,i_s}_{a_s}$ be its $a$-sequence of residues where we assume that for all $j=1,...,s-1$, we have $i_{j}\neq{i_{j+1}}$. Then, we have:

$$A(\underline{\lambda}):=f^{(a_s)}_{i_s}f^{(a_{s-1})}_{i_{s-1}} ...f^{(a_1)}_{i_1}\underline{\emptyset}=\underline{\lambda}+ \sum_{a(\underline{\mu})>a(\underline{\lambda})}{c_{\underline{\lambda},\underline{\mu}}(q)\underline{\mu}},$$
 where $c_{\underline{\lambda},\underline{\mu}}(q)\in{\mathbb{Z}[q,q^{-1}]}$.
\end{prop}
It is obvious that the set $\{A(\underline{\lambda}\ |\ \underline{\lambda}\in{\Lambda}^{1,n},\  n\in{\mathbb{N}}\}$ is a basis of $\overline{\mathcal{M}}_{\mathcal{A}}$. Using the characterization of the canonical basis, we obtain the following theorem:
\begin{thm}\label{base} (\cite[Proposition 4.16]{cyclo}) Let $n\in{\mathbb{N}}$ and let $\underline{\lambda}\in{\Lambda}^{1,n}$, then we have:
$$G({\underline{\lambda}}):=\underline{\lambda}+\sum_{a(\underline{\mu})>a(\underline{\lambda})}b_{\underline{\lambda},\underline{\mu}}(q)\underline{\mu}$$
\end{thm}
In the following paragraph, we provide  an algorithm  which allows us to compute these canonical basis elements.
\section{The algorithm}
We fix $n\in{\mathbb{N}}$, $d\in{\mathbb{N}_{>0}}$, $e\in{\mathbb{N}_{>0}}$ and integers $0\leq v_0\leq v_1\leq ... \leq v_{d-1} <e$. The aim of the algorithm is to compute the decomposition matrix of $\mathcal{H}_{R,n}$ following the proof of \cite[Proposition 4.16]{cyclo}. \\
\\
\textbf{Step 1}: For each $\underline{\lambda}\in{\Lambda^{1,n}}$, we construct the $a$-sequence of residues following Proposition \ref{asuite}:
$$a\textrm{-sequence}(\underline{\lambda})=\underbrace{i_1,...,i_1}_{a_1},\underbrace{i_{2},...,i_{2}}_{a_{2}},...,\underbrace{i_s,...,i_s}_{a_s}.$$
Then, we compute the elements $A(\underline{\lambda})$ of Proposition \ref{decompo} using the  action of Theorem \ref{jmmo}:
$$A(\underline{\lambda}):=f^{(a_s)}_{i_s}f^{(a_{s-1})}_{i_{s-1}} ...f^{(a_1)}_{i_1}\underline{\emptyset}.$$
We obtain a basis   $\{A(\underline{\lambda})\ |\ \underline{\lambda}\in{\Lambda}^1\}$  of $\overline{\mathcal{M}}_{\mathcal{A}}$ which have a ``triangular decomposition''. Since $\overline{f_i}=f_i$ for all $i=0,...,e-1$, we have $\overline{A(\underline{\lambda})}=A(\underline{\lambda})$.\\
\\
\textbf{Step 2}: For each $\underline{\mu}\in{\Pi_n^d}$, we compute its $a$-value\footnote{Note that it is in fact sufficient to compute the values $a_1 (\underline{\mu})$ since we have $a (\underline{\lambda})<a(\underline{\mu})\iff a_1(\underline{\lambda})<a_1(\underline{\mu})$} $a(\underline{\lambda})$ following the definition \ref{afonction}. Let $\underline{\nu}$ be one of the maximal FLOTW $d$-partition with respect to the $a$-function. Then,  we have:
$$G({\underline{\nu}})=A(\underline{\nu}).$$
\textbf{Step 3}: Let  $\underline{\lambda}\in{\Lambda^{1,n}}$. The elements $G({\underline{\mu}})$ with $a(\underline{\mu})>a(\underline{\lambda})$ are known by induction. By   Theorem \ref{base},  there exist polynomials $\alpha_{\underline{\mu},\underline{\lambda}}(q)$ such that:
$$G({\underline{\lambda}})=A(\underline{\lambda})-\sum_{a(\underline{\mu})>a(\underline{\lambda})} \alpha_{\underline{\lambda},\underline{\mu}}(q)G({\underline{\mu}}),\ \  \ \ \ \ \ \ \ \ \ (1)$$
We want to compute  $\alpha_{\underline{\lambda},\underline{\mu}}(q)$ for all $\underline{\mu}\in{\Lambda^{1,n}}$. By Proposition \ref{decompo}, we have:
$$A(\underline{\lambda})= \underline{\lambda}+ \sum_{a(\underline{\mu})>a(\underline{\lambda})}{c_{\underline{\lambda},\underline{\mu}}(q)\underline{\mu}}.$$
Now, since $\overline{G({\underline{\nu}})}=G({\underline{\nu}})$ and $\overline{A(\underline{\nu})}=A(\underline{\nu})$ for all $\underline{\nu}\in{\Lambda^{1,n}}$, we must have  $\alpha_{\underline{\lambda},\underline{\mu}}(q)= \alpha_{\underline{\lambda},\underline{\mu}}(q^{-1})$ for all $\underline{\mu}$  in ${\Lambda^{1,n}}$.

Let $\underline{\nu}\in{\Lambda^{1,n}}$ be one of the minimal  $d$-partition with respect to the $a$-function such that $c_{\underline{\lambda},\underline{\nu}}(q)\notin{q\mathbb{Z}[q]}$. If  $\underline{\nu}$ doesn't exist then, by unicity,  we have $G({\underline{\lambda}})=A(\underline{\lambda})$. If otherwise, by existence of the canonical basis, we have  $\underline{\nu}\in{\Lambda}^{1,n}$. Assume now that we have:
$$c_{\underline{\lambda},\underline{\nu}}(q)=a_{i}q^{i}+a_{i-1}q^{i-1}+...+a_0+...+a_{-i}q^{-i}$$
Where $(a_i)_{j\in{[-i,i]}}$ is a sequence of elements in $\mathbb{Z}$ and $i$ is  a positive integer. Then, we define:
$$\alpha_{\underline{\lambda},\underline{\nu}}(q)=a_{-i}q^{i}+a_{-i+1}q^{i-1}+...+a_0+...+a_{-i}q^{-i}$$
 We have  $\alpha_{\underline{\lambda},\underline{\nu}}(q^{-1})=\alpha_{\underline{\lambda},\underline{\nu}}(q)$. Then, in $(1)$, we replace  $A(\underline{\lambda})$ by  $A(\underline{\lambda})-\alpha_{\underline{\lambda},\underline{\nu}}(q)G({\underline{\nu}})$ which is bar invariant and we repeat this step until  $G({\underline{\lambda}})=A(\underline{\lambda})$.

We finally obtain elements which verify Theorem \ref{basedef} that is to say the canonical basis elements.\\
\\
\textbf{Step 4}:  We  specialize the indeterminate $q$ into $1$ in the canonical basis elements to  obtain the columns of the decomposition matrix of $\mathcal{H}_{R,n}$ which correpond to the   indecomposable  projective $\mathcal{H}_{R,n}$-modules.  
\vspace{0.3cm}

We finally note that we have implemented this algorithm in \textsf{GAP}.

\end{document}